\documentclass[11pt,intlimits]{article}
\usepackage[centertags]{amsmath}
\usepackage[english,german,french]{babel}
\usepackage{amsfonts}
\usepackage{dsfont}
\usepackage{amssymb,amssymb,amscd}
\usepackage{amsthm}
\usepackage[T1]{fontenc}
\usepackage{mathbbol} 
\usepackage{a4wide}
\usepackage{mathrsfs}

\usepackage{color}
\usepackage{stackrel}

\usepackage{multicol}

\usepackage{tikz-cd}
\usepackage{tikz,lipsum}
\usetikzlibrary{shapes}
\usetikzlibrary{decorations.pathmorphing}
\usetikzlibrary{arrows}
\usetikzlibrary{decorations.markings}
\usetikzlibrary{backgrounds,calc}
\tikzset{middlearrow/.style={
        decoration={markings,
            mark= at position 0.6 with {\arrow{#1}} ,
        },
        postaction={decorate}
    }
}
\usepackage{capt-of}
\usepackage{environ}
\makeatletter
\newsavebox{\measure@tikzpicture}
\NewEnviron{scaletikzpicturetowidth}[1]{t%
  \def\tikz@width{#1}%
  \begin{lrbox}{\measure@tikzpicture}%
  \BODY
  \end{lrbox}%
  \pgfmathparse{#1/\wd\measure@tikzpicture}%
  \BODY
}
\makeatother

\usepackage{dynkin-diagrams}
\makeatletter
\newcommand{\extraNode}[6]%
{%
\dynkinPlaceRootRelativeTo{#1}{#2}{#3}{#4}{#5}
\dynkinDefiniteSingleEdge{#1}{#2}
\dynkinRootMark{o}{#1}
\advance\dynkin@nodes by 1
\dynkinLabelRoot{#1}{#6} 
}%
\newcommand{\extraDotNode}[6]%
{%
\dynkinPlaceRootRelativeTo{#1}{#2}{#3}{#4}{#5}
\dynkinIndefiniteSingleEdge{#1}{#2}
\dynkinRootMark{o}{#1}
\advance\dynkin@nodes by 1
\dynkinLabelRoot{#1}{#6} 
}%
\makeatother
\tikzset{/Dynkin diagram,mark=o,edge length=.5cm}

\graphicspath{{Fig/}}

\theoremstyle{thmstyletwo}%
\theoremstyle{thmstylethree}%

\newtheorem{prop}{Proposition}[section]

\newtheorem{rem}[prop]{Remark}
\newtheorem{lem}[prop]{Lemma}
\newtheorem{thm}[prop]{Theorem}
\newtheorem{cor}[prop]{Corollary}

\usepackage{comment}
\usepackage{tikz-cd}
\usepackage{tikz}
\usetikzlibrary{graphs,quotes}
\usepackage[all]{xy}

\begin{document}

\selectlanguage{english}

\title{Existence of a new family of irreducible components in the tensor product and its applications}

\author{Rekha Biswal\footnote{School of Mathematical Sciences, National Institute of Science Education and Research, Bhubaneswar, HBNI, P. O. Jatni, Khurda, 752050, Odisha, India, rekhabiswal27@gmail.com, rekha@niser.ac.in} \ and St\'ephane Gaussent\footnote{Institut Camille Jordan UMR 5208, Université Jean Monnet, CNRS, Centrale Lyon, INSA Lyon, Université Claude Bernard Lyon 1, 20, rue Annino, 42023, Saint-\'Etienne, France, Stephane.gaussent@univ-st-etienne.fr }}

\date{}

\maketitle

\begin{abstract} 

\end{abstract}

\tableofcontents

\abstract{In this paper, using crystal theory, we establish the existence of a new family of irreducible components arising in the tensor product of two irreducible integrable highest weight modules over symmetrizable Kac–Moody algebras. This work is motivated by the Schur positivity conjecture, Kostant’s conjecture, and Wahl’s conjecture. Furthermore, we prove the Schur positivity conjecture in full generality for finite-dimensional simple Lie algebras under the assumption that 
$\lambda  >> \mu$; that is, if $\lambda$ and $\mu$ are the dominant weights in the tensor product, then $\lambda+w\mu$ remains dominant for all $w$ in the Weyl group.

}

\maketitle

\section{Introduction}

Let $\mathfrak{g}$ be a complex symmetrizable Kac–Moody algebra. The irreducible integrable highest weight 
$\mathfrak{g}$-modules are classified by dominant weights. For each dominant weight $\lambda$, we denote by $V(\lambda)$ the corresponding irreducible integrable highest weight module.

For a pair of dominant weights $\lambda$ and $\mu$ of a symmetrizable Kac–Moody algebra $\mathfrak{g}$, a fundamental problem is to determine the decomposition of the tensor product of the corresponding irreducible integrable highest weight modules i.e. to determine the value of $c_{\lambda,\mu}^{\nu}$ in 
$$
V(\lambda)\otimes V(\mu)=\bigoplus_{\nu \in P^{+} } V(\nu)^{\oplus c_{\lambda,\mu}^{\nu}},
$$ where $P^+$ stands for the set of dominant weights.

A wide range of combinatorial, algebraic, and geometric methods are available in the literature for computing tensor product multiplicities
$c_{\lambda,\mu}^{\nu}$. However, even to figure out for which $\nu$, $c_{\lambda,\mu}^{\nu} >0$ using those formulae can be a challenging task. 
In a series of works, various results have been established showing that for certain families of $\nu$, the tensor product multiplicities 
$c_{\lambda,\mu}^{\nu}$
  are strictly positive. For example, the existence of PRV components was proved in \cite{Mathieu}, \cite{Kumar2}, \cite{Kumar3}, \cite{Ressayre}, \cite{Mon} and \cite{MPR}, while the existence of root components was established in \cite{JK} for affine Lie algebras and in \cite{Kumar1} for finite-dimensional simple Lie algebras.

In this paper, using crystal theory, we identify a new family of dominant weights $\nu$ for which $c_{\lambda,\mu}^{\nu} > 0$, and, by applying our main result Theorem \ref{thmSimplyLaced}, we provide affirmative answers to the following three questions in several cases.

\textbf{Question 1 (Schur positivity):} \cite{CFS}
Let $\lambda_1, \lambda_2, \lambda_3, \lambda_4$ be dominant weights for a symmetrizable Kac-Moody algebra $\mathfrak{g}$ such that 
$$
(*)
\left\{
\begin{array}{ll}
\lambda_1 + \lambda_2 = \lambda = \lambda_3 + \lambda_4  \\
\forall\alpha >0, \vert(\lambda_1-\lambda_2)(\alpha^\vee)\vert \leq \vert(\lambda_3 - \lambda_4)(\alpha^\vee)\vert .\\
\end{array}
\right.
$$

Is there an injective \(\mathfrak{g}\)-module homomorphism
\[
V(\lambda_3) \otimes V(\lambda_4) \hookrightarrow V(\lambda_1) \otimes V(\lambda_2)?
\]

This question concerns the comparison of tensor product decompositions via the dominance of differences in weights. A positive answer would imply a natural ordering in the structure of tensor products of irreducible highest weight modules and offer insight into their inclusions within the representation theory of symmetrizable Kac--Moody algebras.

\textbf{Question 2 (Kostant):} \cite{CKM} Let $\rho$ denote the sum of fundamental weights of a symmetrizable Kac-Moody algebra.  Is it true that an irreducible module $V(\beta)$ appears as a component of the tensor product  $V(\rho) \otimes V(\rho)$ if and only if $\beta \leq 2\rho$ in the usual dominance order?

\textbf{Question 3 (Wahl):} \cite{JK}, \cite{Wahl}
Let \(\mathfrak{g}\) be a symmetrizable Kac--Moody algebra with set of positive roots \(R^+\). For dominant weights \(\lambda, \mu \in P^+\), fix an integer \(N \geq 1\). For any dominant weight \(\eta\), define  
\[
S_\eta = \{ i \in I : \langle \eta, \alpha_i^\vee \rangle < N \},
\]  
and for any positive root \(\beta \in R^+\), set  
\[
F_\beta = \{ i \in I : \beta - \alpha_i \notin R^+ \cup \{0\} \}.
\]  
Suppose \(\lambda, \mu \in P^+\) and \(\beta \in R^+\) satisfy the following conditions:  
\begin{enumerate}
    \item \(S_\lambda \cup S_\mu \subseteq F_\beta\),
    \item \(\lambda + \mu - N\beta\) is dominant.
\end{enumerate}  
Is it true that the irreducible module \(V(\lambda + \mu - N\beta)\) appears as a component in the tensor product \(V(\lambda) \otimes V(\mu)\)?

\medskip
The above questions have been extensively studied in various special cases by many authors, yet they remain unresolved in full generality. In particular, the Schur positivity problem continues to be open even for the classical simple Lie algebra $\mathfrak{sl}_n$. The conjecture was originally proposed in \cite{LPP} for $\mathfrak{sl}_n$, and a weaker implication under the Schur positivity hypothesis, namely
\[
c_{\lambda_3,\lambda_4}^\nu > 0 \implies c_{\lambda_1,\lambda_2}^\nu > 0,
\]
was established in \cite{DP}. Subsequently, in \cite{CFS}, the conjecture was extended to arbitrary simple Lie algebras $\mathfrak{g}$ and proved in the special case where the four dominant weights appearing in the hypothesis $(*)$ are scalar multiples of a fixed fundamental weight, as well as for $\mathfrak{g}=\mathfrak{sl}_3$ via crystal-theoretic methods. However, this approach does not easily generalize to broader settings.

Beyond these, the conjecture has been confirmed in several additional special cases by a variety of authors \cite{Fou}, \cite{KL}, \cite{BK}, \cite{FH}, \cite{Naoi1}, \cite{Ven}. These works employ techniques different from crystal theory, notably relying on current algebra methods to explicitly construct generators and relations for the fusion product of irreducible modules. Using these constructions, the authors demonstrate the existence of surjective maps between fusion product modules consistent with the Schur positivity conjecture, thus providing an alternative algebraic viewpoint.

In this paper, we adopt crystal theory to prove the conjecture under a natural and significant additional assumption: that the pair of dominant weights $(\lambda,\mu)$ satisfies the strong dominance condition $\lambda \gg \mu$. Concretely, this means that for every element $w$ in the Weyl group of $\mathfrak{g}$, the weight $\lambda + w\mu$ remains dominant. This condition, while restrictive, is meaningful and has also been considered by Boysal in \cite{Boysal21}, where fusion product multiplicities for $V(\lambda) \otimes V(\mu)$ are computed under similar hypotheses. It remains an intriguing open problem whether the results and techniques developed here, combined with Boysal’s formulae, could be leveraged to establish general inequalities between fusion product multiplicities under the hypothesis $(*)$ of the Schur positivity conjecture.

Additionally, we highlight the work of Di Trani \cite{DiTrani}, who investigated the Kostant conjecture for types $B$, $C$, and $D$, identifying specific irreducible components in the tensor product $V(\rho) \otimes V(\rho)$, thereby contributing to the broader understanding of tensor product decompositions in these classical types.

Beyond the case of strong dominance, we also prove several new results for arbitrary dominant weights $\lambda$ and $\mu$ in the framework of symmetrizable Kac–Moody algebras, as applications of our main theorem. These results provide further supporting evidence toward the validity of the Schur positivity conjecture in more general contexts. A key advantage of our crystal-theoretic approach is its constructive nature: it enables us to explicitly identify some irreducible components appearing in the tensor product with nonzero multiplicity. This contrasts with the fusion product approach, which, while powerful, is more abstract and less explicit in this regard. We believe this explicitness is of independent interest and may lead to further progress in the representation theory of Kac–Moody algebras.

In Section \ref{seExistence}, we establish our main result, Theorem \ref{thmSimplyLaced}, whose proof hinges on a novel formula for certain values of the map $\epsilon$ within crystal theory. We then explore several applications related to the questions introduced above. In Section \ref{se:Deep}, we present a detailed argument proving the Schur positivity conjecture under the assumption $\lambda \gg \mu$.

\medskip

\section{Existence of components in the tensor product}
\label{seExistence}

Let $\mathfrak g$ be a complex symmetrizable Kac-Moody Lie algebra with Cartan subalgebra $\mathfrak h$. Let $A = (a_{i,j})_{i,j\in I}$ be the Cartan matrix of $\mathfrak g$, where $I$ is an index set of cardinality $n$. The set of simple roots of $\mathfrak g$ is denoted by $\{\alpha_i\}_{i\in I}$ and the associated set of simple coroots by $\{\alpha_i^\vee\}_{i\in I}$. Let $Q$ be the root lattice and $P$ the weight lattice of $\mathfrak g$, with $P^{+}$ denoting the set of dominant integral weights.

To each $\lambda \in P^{+}$ corresponds an irreducible integrable highest weight $\mathfrak g$-module $V(\lambda)$.

We consider the crystals $B(\lambda)$ associated to $V(\lambda)$, for each dominant weight $\lambda$. These are semi-normal crystals introduced by Kashiwara \cite{Kashiwara95}. Roughly speaking, $B(\lambda)$ is a set equipped with maps
\[
\mathrm{wt} : B(\lambda) \to P, \quad e_i : B(\lambda) \to B(\lambda) \sqcup \{0\}, \quad f_i : B(\lambda) \to B(\lambda) \sqcup \{0\}, \quad \text{for } i \in I.
\]
The maps $e_i$ and $f_i$ move an element $\pi$ upwards and downwards, respectively, along the string of color $i$ to which $\pi$ belongs. The value $0$ is assigned to $e_i(\pi)$ or $f_i(\pi)$ when $\pi$ is at the upper or lower end of such a string.

The position of $\pi$ in its string of color $i$ is recorded by the functions $\epsilon_i$ and $\phi_i$ defined as
\[
\epsilon_i(\pi) = \max\{n \geq 0 \mid e_i^n(\pi) \neq 0\}, \qquad \phi_i(\pi) = \max\{n \geq 0 \mid f_i^n(\pi) \neq 0\}.
\]
These satisfy compatibility conditions such as
\[
\mathrm{wt}(e_i(\pi)) = \mathrm{wt}(\pi) + \alpha_i, \quad \text{and} \quad \phi_i(\pi) - \epsilon_i(\pi) = \langle \mathrm{wt}(\pi), \alpha_i^\vee \rangle,
\]
for all $\pi \in B(\lambda)$ and $i \in I$.

The tensor product of two crystals $B_1$ and $B_2$,
$B_1 \otimes B_2$, is defined as the Cartesian product of the sets equipped with maps given by
\[
e_i(\pi \otimes \eta) = 
\begin{cases}
e_i(\pi) \otimes \eta & \text{if } \phi_i(\pi) \geq \epsilon_i(\eta), \\
\pi \otimes e_i(\eta) & \text{if } \phi_i(\pi) < \epsilon_i(\eta),
\end{cases}
\]
with a similar definition for $f_i$.

For a given $\lambda \in P^{+}$, we denote the highest weight element in the crystal $B(\lambda)$ by $\pi_\lambda$.

If $\lambda, \mu \in P^{+}$, then it is known that
\[
B(\lambda) \otimes B(\mu) = \bigoplus_{\nu \in P^+} B(\nu)^{\oplus c_{\lambda,\mu}^\nu},
\]
where $c_{\lambda,\mu}^\nu$ is the multiplicity of $V(\nu)$ inside $V(\lambda) \otimes V(\mu)$. Moreover, $c_{\lambda,\mu}^\nu$ equals the number of $\lambda$-dominant elements $\pi$ in $B(\mu)$ of weight $\nu-\lambda$, i.e. \cite{Littelmann95, Littelmann952}
\[
c_{\lambda,\mu}^\nu = \left| \left\{ \pi \in B(\mu) : e_i^{\lambda(\alpha_i^\vee) + 1}(\pi) = 0 \text{ for all } i \in I \text{ and } \mathrm{wt}(\pi) = \nu - \lambda \right\} \right|.
\]

If $\mathfrak g$ is simply laced, Stembridge has defined local axioms for crystals $B(\lambda)$ \cite{S}. We will use two of these axioms, stated as in \cite{BS}:

\begin{itemize}
    \item \textbf{Axiom S1:} For distinct $i,j \in I$, if $x,y \in B(\lambda)$ satisfy $y = e_i x$, then $\epsilon_j(y)$ is either $\epsilon_j(x)$ or $\epsilon_j(x) + 1$. The latter case occurs only if the simple roots $\alpha_i$ and $\alpha_j$ are not orthogonal, i.e., $a_{i,j} \neq 0$.
    \item \textbf{Axiom S2:} For distinct $i,j \in I$, if $x \in B(\lambda)$ satisfies $\epsilon_i(x) > 0$ and $\epsilon_j(e_i x) = \epsilon_j(x) > 0$, then
    \[
    e_i e_j x = e_j e_i x, \quad \text{and} \quad \phi_i(e_j x) = \phi_i(x).
    \]
\end{itemize}

%\medskip
%\begin{lem}\label{lemAnyG}
%Let $\lambda$ be an integral dominant weight of $\mathfrak{g}$. Then
%$$
%f_{i_1}^{b_1}\cdots f_{i_t}^{b_t}\pi_\lambda\neq 0 \iff \lambda(\alpha_{i_k}^\vee) \geq b_k+\sum_{s > k} b_s a_{i_k,i_s}\ ,\ \text{ for all } 1 \leq k \leq t,
%$$  where $i_p \neq i_q$ for $p\neq q$.
%\end{lem}
%\begin{proof}
%Using crystal theory, we see recursively that $f_{i_1}^{b_1}\cdots f_{i_t}^{b_t}\pi_\lambda\neq 0$ if, and only if, $b_k \leq \phi_{i_k}(f_{i_{k+1}}^{b_{k+1}}\cdots f_{i_t}^{b_t}\pi_{\lambda})$ for all $1 \leq k\leq t$. Hence, it is enough to show that $b_k \leq \phi_{i_k}(f_{i_{k+1}}^{b_{k+1}}\cdots f_{i_t}^{b_t}\pi_{\lambda})$ if, and only if, $\lambda(\alpha_{i_k}^\vee) \geq b_k+\sum_{s > k} b_s a_{i_k,i_s}$ for all $1 \leq k\leq t$.  In other words, it is enough to show that $\phi_{i_k}(f_{i_{k+1}}^{b_{k+1}}\cdots f_{i_t}^{b_t}\pi_{\lambda})=\lambda(\alpha_{i_k}^\vee)-\sum_{s>k}b_s a_{i_k,i_s}$. Recall that for any $x \in B(\lambda)$, $\phi_i(x)-\epsilon_i(x)=\text{ wt}(x)(\alpha_i^\vee)$ for all $i\in I$. For $x = f_{i_{k+1}}^{b_{k+1}}\cdots f_{i_t}^{b_t}\pi_{\lambda}$, because $i_p \neq i_q$ for $p\neq q$, one has $\epsilon_{i_k}(f_{i_{k+1}}^{b_{k+1}}\cdots f_{i_t}^{b_t}\pi_{\lambda})=0$. Further, $\text{ wt}(x)(\alpha_{i_k}^\vee) = \lambda(\alpha_{i_k}^\vee) - \sum_{k+1\leq l\leq t} b_l\alpha_{i_l}(\alpha_{i_k}^\vee) =\lambda(\alpha_{i_k}^{\vee})- \sum_{k+1\leq l\leq t} b_l a_{i_k,i_l}$. 

%\end{proof}
%\midskip
\begin{lem}\label{lemAnyG}
Let \(\lambda\) be an integral dominant weight of \(\mathfrak{g}\). Then
\[
f_{i_1}^{b_1} \cdots f_{i_t}^{b_t} \pi_\lambda \neq 0 \iff \lambda(\alpha_{i_k}^\vee) \geq b_k + \sum_{s > k} b_s a_{i_k, i_s} \quad \text{for all } 1 \leq k \leq t,
\]
where \(i_p \neq i_q\) for all \(p \neq q\).
\end{lem}

\begin{proof}
Using crystal theory, we proceed recursively. The element
\[
f_{i_1}^{b_1} \cdots f_{i_t}^{b_t} \pi_\lambda \neq 0
\]
if and only if
\[
b_k \leq \varphi_{i_k}(f_{i_{k+1}}^{b_{k+1}} \cdots f_{i_t}^{b_t} \pi_\lambda)
\quad \text{for all } 1 \leq k \leq t.
\]
Hence, it suffices to show that for each \(k\),
\[
b_k \leq \varphi_{i_k}(f_{i_{k+1}}^{b_{k+1}} \cdots f_{i_t}^{b_t} \pi_\lambda) \iff \lambda(\alpha_{i_k}^\vee) \geq b_k + \sum_{s > k} b_s a_{i_k, i_s}.
\]

In other words, it suffices to prove that
\[
\varphi_{i_k}(f_{i_{k+1}}^{b_{k+1}} \cdots f_{i_t}^{b_t} \pi_\lambda) = \lambda(\alpha_{i_k}^\vee) - \sum_{s > k} b_s a_{i_k, i_s}.
\]

Recall that for any \(x \in B(\lambda)\), we have the relation
\[
\varphi_i(x) - \epsilon_i(x) = \mathrm{wt}(x)(\alpha_i^\vee) \quad \text{for all } i \in I.
\]

For
\[
x = f_{i_{k+1}}^{b_{k+1}} \cdots f_{i_t}^{b_t} \pi_\lambda,
\]
and since the indices \(i_p\) are distinct, we have
\[
\epsilon_{i_k}(x) = 0.
\]

Furthermore, the weight of \(x\) evaluated at \(\alpha_{i_k}^\vee\) is
\[
\mathrm{wt}(x)(\alpha_{i_k}^\vee) = \lambda(\alpha_{i_k}^\vee) - \sum_{l = k+1}^t b_l \alpha_{i_l}(\alpha_{i_k}^\vee) = \lambda(\alpha_{i_k}^\vee) - \sum_{l = k+1}^t b_l a_{i_k, i_l}.
\]

Combining these, we conclude
\[
\varphi_{i_k}(x) = \mathrm{wt}(x)(\alpha_{i_k}^\vee) + \epsilon_{i_k}(x) = \lambda(\alpha_{i_k}^\vee) - \sum_{l > k} b_l a_{i_k, i_l}.
\]

This completes the proof.
\end{proof}

%\midskip
\subsection{Simply laced cases}
In this subsection, let $\mathfrak{g}$ be a complex simply-laced Kac-Moody Lie algebra of rank $n$. We denote by $I = \{1, 2, \ldots, n\}$ the indexing set of cardinality $n$ corresponding to the nodes of the Dynkin diagram. The entries $a_{i,j}$ of the Cartan matrix of $\mathfrak{g}$, for $i,j \in I$, take values in $\{0, -1, 2\}$. The following theorem is the main result of this paper.

\begin{thm}\label{thmSimplyLaced}
Let $t$ be a positive integer, and let $\lambda, \mu$ be dominant integral weights of $\mathfrak{g}$. Suppose there exist $(b_1, \ldots, b_t) \in \mathbb{N}^t$ and $(i_1, \ldots, i_t) \in I^t$ satisfying the conditions:
\[
\lambda(\alpha_{i_t}^\vee) \geq b_t, \quad \mu(\alpha_{i_1}^\vee) \geq b_1,
\]
\[
\mu(\alpha_{i_r}^\vee) \geq b_r + \sum_{s<r} b_s a_{i_s, i_r}, \quad \text{for all } 1 < r \leq t,
\]
\[
\lambda(\alpha_{i_r}^\vee) \geq b_r + \sum_{s>r} b_s a_{i_s, i_r}, \quad \text{for all } 1 \leq r < t.
\]
Then, provided that $p \neq q$ implies $\alpha_{i_p} \neq \alpha_{i_q}$, the module 
\[
V\Big(\lambda + \mu - \sum_{p=1}^t b_p \alpha_{i_p}\Big)
\]
is an irreducible component of the tensor product $V(\lambda) \otimes V(\mu)$.
\end{thm}

The proof of Theorem \ref{thmSimplyLaced} is split into several intermediate results.

\begin{prop}\label{LemEpsilon}
Let $\lambda$ be a dominant weight of $\mathfrak{g}$. Let $(i_1,\ldots, i_t) \in I^t$ be such that $p \neq q$ implies $i_p \neq i_q$. Then for $1 \leq r \leq t$ and for $j \notin \{i_1,\ldots,i_t\}$, provided $f_{i_1}^{b_1} \cdots f_{i_t}^{b_t} \pi_\lambda \neq 0$, we have
\[
\epsilon_{i_r}(f_{i_1}^{b_1} \cdots f_{i_t}^{b_t} \pi_\lambda) = \max \left\{ 0, b_r + \sum_{s < r} b_s a_{i_s,i_r} \right\}, \quad \epsilon_j(f_{i_1}^{b_1} \cdots f_{i_t}^{b_t} \pi_\lambda) = 0.
\]
\end{prop}

\begin{proof}
We prove this by induction on $m = \sum_{k=1}^t b_k$. The statement holds if $m=1$.

Assume it is true for $m$, and consider the case $m+1$. The statement is clear for $r=1$. Let $r > 1$ and set
\[
x = f_{i_1}^{b_1} \cdots f_{i_t}^{b_t} \pi_\lambda, \quad y = e_{i_1} x = f_{i_1}^{b_1 - 1} \cdots f_{i_t}^{b_t} \pi_\lambda.
\]
By the induction hypothesis,
\[
\epsilon_{i_r}(y) = \max \left\{ 0, b_r + \sum_{1 < s < r} b_s a_{i_s,i_r} + (b_1 - 1) a_{i_1,i_r} \right\}.
\]

\medskip
\textbf{CASE 1:} Assume that $\epsilon_{i_r}(y) > 0$.

If $a_{i_1,i_r} = 0$, then by Axiom S1, $\epsilon_{i_r}(y) = \epsilon_{i_r}(x)$, which proves the statement.

If $a_{i_1,i_r} = -1$, then $\epsilon_{i_r}(y) = \epsilon_{i_r}(x) + 1$. Indeed, if $\epsilon_{i_r}(y) = \epsilon_{i_r}(x) > 0$, then by Axiom S2,
\[
\phi_{i_1}(e_{i_r} x) = \phi_{i_1}(x).
\]
Recall that
\[
\phi_{i_1}(x) = \epsilon_{i_1}(x) + \operatorname{wt}(x)(\alpha_{i_1}^\vee) = b_1 + \operatorname{wt}(x)(\alpha_{i_1}^\vee),
\]
and
\[
\phi_{i_1}(e_{i_r} x) = \epsilon_{i_1}(e_{i_r} x) + \operatorname{wt}(x)(\alpha_{i_1}^\vee) - 1,
\]
implying $\epsilon_{i_1}(e_{i_r} x) = b_1 + 1$.

This is impossible because 
\[
\operatorname{wt}(e_{i_r} x) = \lambda - \sum_{u=1}^t b_u \alpha_{i_u} + \alpha_{i_r}
\]
implies 
\[
\epsilon_{i_1}(e_{i_r} x) \leq b_1,
\]
since otherwise, 
\[
\operatorname{wt}(e_{i_1}^{b_1 + 1} e_{i_r} x) = \lambda - \sum_{1 < u \leq t} b_u \alpha_{i_u} + \alpha_{i_1} + \alpha_{i_r}
\]
would not be smaller than $\lambda$.

\medskip
\textbf{CASE 2:} Assume that $\epsilon_{i_r}(y) = 0$.

By Proposition 4.5 of \cite{BS}, if $\epsilon_{i_1}(x) > 0$, then
\[
\epsilon_{i_r}(x) = \epsilon_{i_r}(y) - 1 \quad \text{or} \quad \epsilon_{i_r}(x) = \epsilon_{i_r}(y).
\]
Since $\epsilon_{i_r}(y) = 0$ and $\epsilon_{i_r}(x) \geq 0$, it follows that $\epsilon_{i_r}(x) = 0$, proving the statement.
\end{proof}

\begin{lem}\label{lemEquivDom}
Let $\lambda$ and $\mu$ be dominant weights of $\mathfrak{g}$. Let $(i_1, \ldots, i_t) \in I^t$ satisfy $i_p \neq i_q$ for $p \neq q$. Then $x = f_{i_1}^{b_1} \cdots f_{i_t}^{b_t} \pi_\lambda \neq 0$ is $\mu$-dominant if and only if the following conditions hold:
\[
\lambda(\alpha_{i_t}^\vee) \geq b_t, \quad \mu(\alpha_{i_1}^\vee) \geq b_1,
\]
\[
\mu(\alpha_{i_r}^\vee) \geq b_r + \sum_{s < r} b_s a_{i_s, i_r}, \quad \text{for all } 1 < r \leq t,
\]
\[
\lambda(\alpha_{i_r}^\vee) \geq b_r + \sum_{s > r} b_s a_{i_s, i_r}, \quad \text{for all } 1 \leq r < t.
\]
\end{lem}

\begin{proof}
The conditions on $\lambda$ follow from Lemma \ref{lemAnyG}.

For the conditions on $\mu$, recall that $x$ is $\mu$-dominant if and only if for all $i \in I$,
\[
e_i^{\mu(\alpha_i^\vee) + 1}(x) = 0.
\]
For $i \notin \{i_1, \ldots, i_t\}$, we have $\epsilon_i(x) = 0$, so the condition is automatic.

For $i = i_r \in \{i_1, \ldots, i_t\}$, Proposition \ref{LemEpsilon} says that
\[
\epsilon_{i_r}(x) = \max \left\{ 0, b_r + \sum_{s < r} b_s a_{i_s, i_r} \right\}.
\]
Hence, $x$ is $\mu$-dominant if and only if
\[
\mu(\alpha_{i_r}^\vee) \geq b_r + \sum_{s < r} b_s a_{i_s, i_r} \quad \text{for all } 1 \leq r \leq t,
\]
which completes the proof.
\end{proof}

These two results prove Theorem \ref{thmSimplyLaced}.

\begin{cor}
Let $\lambda$ and $\mu$ be two regular dominant weights of a simply laced symmetrizable Kac-Moody algebra. Let $\alpha$ be an element in the root lattice of $\mathfrak{g}$ which is a sum of distinct simple roots. Then $V(\lambda + \mu - \alpha)$ is an irreducible component of $V(\lambda) \otimes V(\mu)$.

In particular, if $\rho$ is the sum of fundamental weights for a simply laced Kac-Moody algebra $\mathfrak{g}$, then for any such $\alpha$, $V(2\rho - \alpha)$ is an irreducible component of $V(\rho) \otimes V(\rho)$.
\end{cor}

\begin{proof}
This follows from the fact that if $\alpha = \alpha_{i_1} + \cdots + \alpha_{i_t}$ with $i_1 < \cdots < i_t$, then the element
\[
f_{i_1} \cdots f_{i_t} \pi_\mu \in B(\mu)
\]
is always nonzero and $\lambda$-dominant.
\end{proof}

%\midskip
\begin{rem}

Note that $2\rho - \alpha$ is a PRV component of $V(\rho) \otimes V(\rho)$ for some roots $\alpha$, but not for all $\alpha$ that are sums of distinct simple roots. For example, if the Dynkin diagram of $\mathfrak{g}$ contains a cyclic subdiagram with an odd number of nodes, and if we take $\alpha$ to be the sum of the simple roots corresponding to those nodes, then $2\rho - \alpha$ can never be expressed in the form $\sigma \rho + \tau \rho$ for any Weyl group elements $\sigma, \tau$ of $\mathfrak{g}$. Hence, this corollary identifies some new components inside $V(\rho) \otimes V(\rho)$, providing additional supporting evidence towards Kostant’s conjecture.
\end{rem}

\subsection{Applications to Schur positivity}

Let us first mention that the hypothesis $(*)$ in the Schur positivity conjecture can also be written as the right hand side:
$$
\left\{
\begin{array}{ll}
\lambda_1 + \lambda_2 = \lambda = \lambda_3 + \lambda_4 \quad \text{and } \forall\alpha >0, \\
\vert(\lambda_1-\lambda_2)(\alpha^\vee)\vert \leq \vert(\lambda_3 - \lambda_4)(\alpha^\vee)\vert \\
\end{array}
\right.
\Leftrightarrow
\left\{
\begin{array}{ll}
\lambda_1 + \lambda_2 = \lambda = \lambda_3 + \lambda_4 \quad \text{and } \forall\alpha >0,  \\
\text{min}\{\lambda_3(\alpha^{\vee}),\lambda_4(\alpha^{\vee})\} \leq \text{min}\{\lambda_1(\alpha^{\vee}),\lambda_2(\alpha^{\vee})\} .\\
\end{array}
\right.
$$

\begin{prop}\label{Schur}

Let $\lambda_1,\lambda_2,\lambda_3,\lambda_4$ be dominant integral weights of $\mathfrak{g}$ that satisfy the hypothesis $(*)$.
Assume that $(i_1,\cdots,i_t) \in I^t$ is such that $r\neq s$ implies $i_r\neq i_s$ and that the subdiagram spanned by the nodes $i_1, ... , i_t$ is connected of type $A,D,E$. Assume also that $f_{i_1}\cdots f_{i_t}\pi_{\lambda_4}$ is $\lambda_3$-dominant. Then there exists a permutation $\sigma\in S_t$ such that $f_{i_{\sigma(1)}}\cdots f_{i_{\sigma(t)}}\pi_{\lambda_2}$ is nonzero and $\lambda_1$-dominant.

\end{prop}

\begin{proof}
It is clear from  Lemma \ref{lemEquivDom} that for any two dominant weights $\lambda$ and $\mu$, if $p_\ell \neq p_m$ for $\ell \neq m$ then the element $f_{p_1}\cdots f_{p_t}\pi_{\lambda}$ is nonzero and $\mu$-dominant iff for all $1 \leq u \leq t$:
\begin{itemize}
\item[(1)] $\mu(\alpha_{p_u}^{\vee}) = 0$ $\Longrightarrow$ there exists $v < u$ such that $a_{p_u,p_v} \neq 0$.
\item[(2)]$\lambda(\alpha_{p_u}^{\vee}) = 0$ $\Longrightarrow$ there exists $w > u$ such that $a_{p_u,p_w} \neq 0$. 
\end{itemize}
In particular, for $1 \leq k \leq t$, if $|\{q: 1 \leq q \leq t \text{ and } a_{p_q,p_k} = -1\}|=1$, then either $\lambda(\alpha_{p_k}^{\vee}) \geq 1$ or $\mu(\alpha_{p_k}^{\vee}) \geq 1$ which implies that $(\lambda+\mu)(\alpha_{p_k}^{\vee}) \geq 1$.

More precisely, for $1 \leq k \leq t$, a node $p_k$ violates the required conditions for $f_{p_1}\cdots f_{p_t}\pi_{\lambda}$ to be nonzero and $\mu$-dominant if and only if either $\mu(\alpha_{p_k}^{\vee})=0$ and there is no connected node to the left of $p_k$ in the monomial $f_{p_1}\cdots f_{p_t}$ (i.e. $a_{p_u,p_k}=0$  for all $u<k$ ) or $\lambda(\alpha_{p_k}^{\vee})=0$ and there is no connected node to the right of $p_k$ in the monomial $f_{p_1}\cdots f_{p_t}$ ( i.e.  $a_{p_v,p_k}=0$ for all $v>k$). We will be using these facts in the proof below.

Let us assume that $t \geq 2$ as the statement is obvious for $t=1$. Using a permutation in $S_t$, we can assume that $i_1, i_2,..., i_t$ are nodes of the Dynkin diagram, labeled as in Kac's book \cite{Kac}.

Assume first that the nodes $i_1,\cdots,i_t$ form a Dynkin subdiagram of type $A$ as in the diagram below such that $\{i_1,\cdots ,i_t\}=\{j_1,\cdots,j_t\}$. 
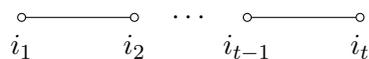
\begin{center}
\begin{tikzpicture}
[edot/.style={fill=white,circle,draw,minimum size=3pt,inner sep=0,outer sep=0},xscale=1.5]  

\path 
(0,0) node{$\cdots$}
(.5,0) node[edot] (R1) {} node[below=1mm]{$j_{t-1}$}
(1.5,0) node[edot] (R2) {}node[below=1mm] {$j_t$}
(-.5,0) node[edot] (L1) {} node[below=1mm] {$j_2$}
(-1.5,0) node[edot] (L2) {} node[below=1mm] {$j_1$}
;       
\draw (L1)--(L2) (R1)--(R2);
\end{tikzpicture}
\end{center}
 Since $|\{q: 1 \leq q \leq t \text{ and } a_{j_q,j_1}=-1\}| = |\{q: 1 \leq q \leq t \text{ and } a_{j_q,j_t}= -1\}|=1$, the conditions in the Lemma \ref{lemEquivDom} for a nonzero element $f_{i_1}\cdots f_{i_t}\pi_{\lambda_4}$ to be $\lambda_3$-dominant implies that 
$(\lambda_3+\lambda_4)(\alpha_{j_1}^{\vee}) \geq 1$ and $(\lambda_3+\lambda_4)(\alpha_{j_t}^{\vee}) \geq 1$. But $\lambda_1+\lambda_2=\lambda_3+\lambda_4$ implies that $(\lambda_1+\lambda_2)(\alpha_{j_1}^{\vee}) \geq 1$ and $(\lambda_1+\lambda_2)(\alpha_{j_t}^{\vee}) \geq 1$ must be true. So we have the following four different cases:

\textbf{Case 1}: Assume that $\lambda_1(\alpha^\vee_{j_1})\geq 1$ and $\lambda_2(\alpha^\vee_{j_t})\geq 1$. In this case, Lemma \ref{lemEquivDom} implies that $f_{j_1}\cdots f_{j_t}\pi_{\lambda_2}$ is nonzero and $\lambda_1$-dominant. This follows from the fact that $a_{j_1,j_2} \neq 0$ and $a_{j_{t-1},j_t} \neq 0$ and in the notation of Lemma \ref{lemEquivDom}, $b_i=1$ for all $1 \leq i \leq t$ for the monomial $f_{j_1}\cdots f_{j_t}\pi_{\lambda_2}$.
Hence for any $1<r<t$, $b_r + \sum_{s<r} b_sa_{j_s,j_r} = 1 +  a_{j_{r-1},j_r} = 1-1 = 0$. Similarly, for any $1<r<t$, $b_r + \sum_{s>r} b_sa_{j_s,j_r} = 1 +  a_{j_{r+1},j_r} = 1-1 = 0$. Therefore the only required conditions for  $f_{j_1}\cdots f_{j_t}\pi_{\lambda_2}$ to be nonzero and $\lambda_1$-dominant in this case is $\lambda_1(\alpha^\vee_{j_r})\geq 0$ and $\lambda_2(\alpha^\vee_{j_r})\geq 0$ for $1 < r < t$. However these two conditions are automatically satisfied because $\lambda_1$ and $\lambda_2$ are dominant weights. 

\textbf{Case 2}: assume that $\lambda_2(\alpha^\vee_{j_1})\geq 1$ and $\lambda_1(\alpha^\vee_{j_t})\geq 1$.\\
This case is analogous to the previous one and in this case it is easy to check that $f_{j_t}\cdots f_{j_1}\pi_{\lambda_2}$ is nonzero and $\lambda_1$-dominant by Lemma \ref{lemEquivDom}.
 
%\medskip
\textbf{Case 3}: assume that $\lambda_1(\alpha^\vee_{j_1})\geq 1$, $\lambda_2(\alpha^\vee_{j_1})=0$ and $\lambda_1(\alpha^\vee_{j_t})\geq 1$, $\lambda_2(\alpha^\vee_{j_t})=0$. \\
Because $f_{i_1}\cdots f_{i_t}\pi_{\lambda_4}$ is nonzero and $\lambda_3$-dominant, we have $\lambda_3(\alpha^\vee_{i_1}) \geq 1$ and $\lambda_4(\alpha^\vee_{i_t}) \geq 1$. So, for $\beta = \alpha_{i_1} + \cdots +\alpha_{i_t}=\alpha_{j_1} + \cdots +\alpha_{j_t}$, min$\{\lambda_3(\beta^\vee),\lambda_4(\beta^\vee)\}\geq 1$. \\
Then, $1 \leq\text{min}\{\lambda_3(\beta^{\vee}),\lambda_4(\beta^{\vee})\} \leq \text{min}\{\lambda_1(\beta^{\vee}),\lambda_2(\beta^{\vee})\}$ for which $\lambda_2(\beta^{\vee}) \geq 1$. Hence, the fact that $\lambda_2(\alpha^\vee_{j_1} + \alpha^\vee_{j_t}) = 0$ implies that there exists a $1<k<t$ such that $\lambda_2(\alpha^\vee_{j_k})\geq 1$. Choose such a $k$ maximal, then $f_{j_1}\cdots f_{j_{k-1}}f_{j_t}\cdots f_{j_k}\pi_{\lambda_2}$ is nonzero and $\lambda_1$-dominant. 

%Indeed, for $2 \\leq r\\\leq k-1$, we have $1 + \sum_{s<r} a_{j_s,j_r} = 1 + a_{j_{r-1},j_r} = 1-1 = 0$. 

Indeed, for $2 \leq r\leq k-1$ and $k+1 \leq r \leq t-1$, there exists one connected node to the left of $j_r$ and another connected node to the right of $j_r$ in the element $f_{j_1}\cdots f_{j_{k-1}}f_{j_t}\cdots f_{j_k}\pi_{\lambda_2}$. Hence $\lambda_1(\alpha^\vee_{j_r})\geq 0$ and $\lambda_2(\alpha^\vee_{j_r})\geq 0$ for such $r$ implies that the conditions of Lemma \ref{lemEquivDom} are satisfied for the nodes $i_r$ when $2 \leq r \leq k-1$ and $k+1 \leq r \leq t-1$. For $j_k$, because $\lambda_2(\alpha_{j_k}^{\vee})\geq 1$ and there is a connected node to the left of $j_k$ in the element $f_{j_1}\cdots f_{j_{k-1}}f_{j_t}\cdots f_{j_k}\pi_{\lambda_2}$, the required conditions are also satisfied for the node $j_k$. Similarly for $j_1, j_t$, because $\lambda_1(\alpha^\vee_{j_1})\geq 1, \lambda_1(\alpha^\vee_{j_t})\geq 1$ and there are connected nodes to the right of both $j_1$ and $j_t$ in the element $f_{j_1}\cdots f_{j_{k-1}}f_{j_t}\cdots f_{j_k}\pi_{\lambda_2}$, the required conditions are also satisfied for the nodes $j_1$ and $j_t$.

%Further, for $r\neq k $, we also have $b_r + \sum_{s>r} b_sa_{j_s,j_r} = 1-1 = 0$. So the conditions of Lemma \ref{lemEquivDom} are satisfied for $\lambda_2$. 

%\medskip
\textbf{Case 4}: assume that $\lambda_1(\alpha^\vee_{j_1})= 0$, $\lambda_2(\alpha^\vee_{j_1})\geq 1$ and $\lambda_1(\alpha^\vee_{j_t}) = 0$, $\lambda_2(\alpha^\vee_{j_t})\geq 1$. \\
This case is analogous to the previous one. Because $f_{i_1}\cdots f_{i_t}\pi_{\lambda_4}$ is nonzero and $\lambda_3$-dominant, we have $\lambda_3(\alpha^\vee_{i_1}) \geq 1$ and $\lambda_4(\alpha^\vee_{i_t}) \geq 1$. 

So, for $\beta = \alpha_{i_1} + \cdots +\alpha_{i_t}=\alpha_{j_1} + \cdots +\alpha_{j_t}$, min$\{\lambda_3(\beta^\vee),\lambda_4(\beta^\vee)\}\geq 1$. Then, $1 \leq\text{min}\{\lambda_3(\beta^{\vee}),\lambda_4(\beta^{\vee})\} \leq \text{min}\{\lambda_1(\beta^{\vee}),\lambda_2(\beta^{\vee})\}$. Hence, the fact that $\lambda_1(\alpha^\vee_{j_1} + \alpha^\vee_{j_t}) = 0$ implies that there exists a $1<k<t$ such that $\lambda_1(\alpha^\vee_{j_k})\geq 1$. Choose such a $k$ minimal, then $f_{j_k}\cdots f_{j_1}f_{j_{k+1}}\cdots f_{j_t}\pi_{\lambda_2}$ is nonzero and $\lambda_1$-dominant.

\medskip
(Note that in both cases 3 and 4, we could choose any $k$ satisfying the required condition and not necessarily minimal or maximal.)
\medskip

The case where $\lambda_1(\alpha^\vee_{j_1})=\lambda_2(\alpha^\vee_{j_1}) = 0$ or $\lambda_1(\alpha^\vee_{j_t}) = \lambda_2(\alpha^\vee_{j_t}) = 0$ can not arise. This is because $|\{q: 1 \leq q \leq t \text{ and } a_{j_q,j_t} = -1\}|=|\{q: 1 \leq q \leq t \text{ and } a_{j_q,j_1} = -1\}|=1$ and $\{i_1, \cdots, i_t\}=\{j_1,\cdots,j_t\}$. Hence $f_{i_1}\cdots f_{i_t}\pi_{\lambda_4}$ is nonzero and $\lambda_3$-dominant implies that $(\lambda_3+\lambda_4)(\alpha_{j_1}^{\vee}) \geq 1 $ and $ (\lambda_3+\lambda_4)(\alpha_{j_t}^{\vee}) \geq 1$ must be true. Since we assume that $\lambda_3+\lambda_4=\lambda_1+\lambda_2$, $(\lambda_1+\lambda_2)(\alpha_{j_1}^{\vee}) \geq 1$ and $(\lambda_1+\lambda_2)(\alpha_{j_t}^{\vee}) \geq 1$ must be true as well.

\medskip
Now, let us consider the case when the Dynkin diagram spanned by the nodes $i_1,\cdots,i_t$ is as follows which includes the case of type D and E.

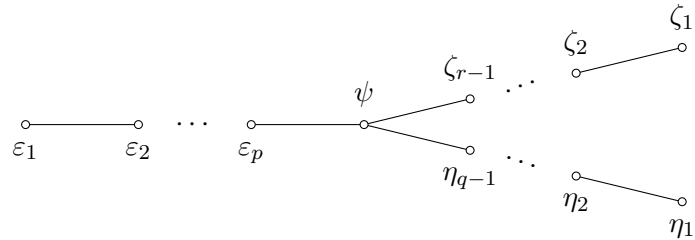
\begin{center}
\begin{tikzpicture}
[edot/.style={fill=white,circle,draw,minimum size=3pt,inner sep=0,outer sep=0},xscale=1.5]  

\def\t{20}
\path 
(0,-4) node[edot] (3) {} node[above=1mm] {$\psi$}
++(\t:1) node[edot] (3R1a) {} node[above=1mm]{$\zeta_{r-1}$}
++(\t:.5) node[rotate=\t]{$\cdots$}
++(\t:.5) node[edot] (3R2a) {} node[above=1mm]{$\zeta_2$}
++(\t:1) node[edot] (3R3a) {} node[above=1mm]{$\zeta_1$}
(0,-4)
++(-\t:1) node[edot] (3R1b) {} node[below=1mm]{$\eta_{q-1}$}
++(-\t:.5) node[rotate=-\t]{$\cdots$}
++(-\t:.5) node[edot] (3R2b) {} node[below=1mm]{$\eta_2$}
++(-\t:1) node[edot] (3R3b) {} node[below=1mm]{$\eta_1$}
(-1,-4) node[edot] (3L1) {} node[below=1mm]{$\varepsilon_{p}$}
(-1.5,-4) node {$\cdots$}
(-2,-4) node[edot] (3L2) {} node[below=1mm]{$\varepsilon_2$}
(-3,-4) node[edot] (3L3) {} node[below=1mm]{$\varepsilon_1$}
;
\draw (3)--(3R1a) (3R2a)--(3R3a) (3)--(3R1b) (3R2b)--(3R3b)
(3)--(3L1) (3L2)--(3L3);
\end{tikzpicture}
\end{center}

The proof proceeds in close analogy to the type~$A$ case. As before, it is clear that $$(\lambda_1+\lambda_2)(\alpha^\vee_{\eta_1}) \geq 1, \ \ (\lambda_1+\lambda_2)(\alpha^\vee_{\zeta_1}) \geq1, \ \
(\lambda_1+\lambda_2)(\alpha^\vee_{\varepsilon_1}) \geq 1.$$ 
Indeed, this follows from the fact that there is only one connected node to each of the nodes $\eta_1,\zeta_1$ and $\varepsilon_1$ in the above Dynkin diagram for which $(\lambda_3+\lambda_4)(\alpha^\vee_{\eta_1}) \geq 1, \ (\lambda_3+\lambda_4)(\alpha^\vee_{\zeta_1}) \geq1, \ 
(\lambda_3+\lambda_4)(\alpha^\vee_{\varepsilon_1}) \geq 1$ must be true. We distinguish cases based on the number of elements in the set  \{$\lambda_1(\alpha_{\varepsilon_1}^{\vee}), \lambda_2(\alpha_{\varepsilon_1}^{\vee}), \lambda_1(\alpha_{\eta_1}^{\vee}),\lambda_2(\alpha_{\eta_1}^{\vee}), \lambda_1(\alpha_{\zeta_1}^{\vee}), \lambda_2(\alpha_{\zeta_1}^{\vee})\}$ that are equal to zero.  We present the proof for certain representative cases, as the remaining cases follow by analogous arguments.

\textbf{Case 1}: $\lambda_1(\alpha_{\varepsilon_1}^{\vee}) \geq 1, \lambda_2(\alpha_{\varepsilon_1}^{\vee})=0$, $\lambda_1(\alpha_{\zeta_1}^{\vee}) \geq 1, \lambda_2(\alpha_{\zeta_1}^{\vee})=0$ ,$\lambda_1(\alpha_{\eta_1}^{\vee}) = 0, \lambda_2(\alpha_{\eta_1}^{\vee})\geq 1$\\
In this case, one can see that $f_{\varepsilon_1}\cdots f_{\varepsilon_p}  f_{\zeta_{1}}\cdots f_{\zeta_{r-1}}f_{\psi}f_{\eta_{q-1}}\cdots f_{\eta_1} \pi_{\lambda_2}$ is nonzero and \\ $\lambda_1$-dominant.
This follows from the fact that except the boundary $f$-operators i.e. $f_{\varepsilon_1}$ and $f_{\eta_1}$ of the monomial $f_{\varepsilon_1}\cdots f_{\varepsilon_p}  f_{\zeta_{1}}\cdots f_{\zeta_{r-1}}f_{\psi}f_{\eta_{q-1}}\cdots f_{\eta_1}$, all the intermediate $f$-operators except $f_{\zeta_1}$ has a connected node to the left and another connected node to the right of it. Again nodes $\varepsilon_1$,$\eta_1$ and $\zeta_1$ do not violate the conditions in this case. That follows from the fact that: \begin{itemize}
\item[(1)] $\lambda_1(\alpha_{\varepsilon_1}^{\vee}) \geq 1, \lambda_2(\alpha_{\varepsilon_1}^{\vee})=0$ and there is a connected node to the right of $\varepsilon_1$
\item[(2)] $\lambda_1(\alpha_{\eta_1}^{\vee}) = 0, \lambda_2(\alpha_{\eta_1}^{\vee})\geq 1$ and there is a connected node to the left of $\eta_1$. 
\item [(3)]$\lambda_1(\alpha_{\zeta_1}^{\vee}) \geq 1, \lambda_2(\alpha_{\zeta_1}^{\vee})=0$ and there is a connected node to the right of $\zeta_1$ .
\end{itemize}

\textbf{Case 2}: $\lambda_1(\alpha_{\varepsilon_1}^{\vee}) =0, \lambda_2(\alpha_{\varepsilon_1}^{\vee})\geq 1$, $\lambda_1(\alpha_{\zeta_1}^{\vee})=0, \lambda_2(\alpha_{\zeta_1}^{\vee})\geq 1$ ,$\lambda_1(\alpha_{\eta_1}^{\vee}) \geq 1, \lambda_2(\alpha_{\eta_1}^{\vee})=0$\\
In this case, $f_{\eta_1}\cdots f_{\eta_{q-1}}f_{\psi} f_{\varepsilon_q}\cdots f_{\varepsilon_1} f_{\zeta_{r-1}}\cdots f_{\zeta_1} \pi_{\lambda_2}$ is nonzero and $\lambda_1$-dominant because of the analogous reasoning.

\textbf{Case 3}: $\lambda_1(\alpha_{\varepsilon_1}^{\vee}) = 0, \lambda_2(\alpha_{\varepsilon_1}^{\vee}) \geq 1$, $\lambda_1(\alpha_{\eta_1}^{\vee}) \geq 1 $, $\lambda_2(\alpha_{\eta_1}^{\vee}) = 0, $  $ \lambda_1(\alpha_{\zeta_1}^{\vee}) \geq 1,\lambda_2(\alpha_{\zeta_1}^{\vee}) = 0 $\\
In this case, similarly one can check that $f_{\zeta_1} \cdots f_{\zeta_{r-1}} f_{\eta_1} \cdots f_{\eta_{q-1}} f_{\psi} f_{\varepsilon_p} \cdots f_{\varepsilon_1} \pi_{\lambda_2}$ is nonzero and $\lambda_1$-dominant.

\textbf{Case 4}: $\lambda_2(\alpha_{\eta_1}^\vee)=\lambda_2(\alpha_{\zeta_1}^\vee)=\lambda_2(\alpha_{\varepsilon_1}^\vee)=0$ and $\lambda_1(\alpha_{\eta_1}^\vee) \geq 1, \lambda_1(\alpha_{\zeta_1}^\vee) \geq 1, \lambda_1(\alpha_{\varepsilon_1}^\vee)\geq 1$\\
As before, $c_{\lambda_3,\lambda_4}^{\lambda_3+\lambda_4-\beta} > 0$, for $\beta=\alpha_{\eta_1}+\cdots+\alpha_{\eta_{q-1}}+\alpha_{\psi}+\alpha_{\zeta_1}+\cdots+\alpha_{\zeta_{r-1}}+\alpha_{\varepsilon_1}+\cdots+\alpha_{\varepsilon_p}$, implies that $\lambda_3(\beta^\vee) \geq 1$ and $\lambda_4(\beta^\vee) \geq 1$ for which $\text{min} \{\lambda_1(\beta^\vee),\lambda_2(\beta^\vee)\} \geq 1$. Hence there must exists $k >1$ such that at least one of $\lambda_2(\alpha_{\eta_k}^\vee) \geq 1$ or $\lambda_2(\alpha_{\zeta_k}^\vee) \geq 1$ or $\lambda_2(\alpha_{\varepsilon_k}^\vee) \geq 1$ or $\lambda_2(\alpha_{\psi}^\vee) \geq 1$ is true. Depending on which holds, the following elements are nonzero and $\lambda_1$-dominant.
\begin{itemize}
\item[(1)] If $\lambda_2(\alpha^\vee_{\zeta_k}) \geq 1$, then $f_{\eta_1}\cdots f_{\eta_{q-1}}f_{\varepsilon_1}\cdots f_{\varepsilon_p} f_{\psi} f_{\zeta_{r-1}} \cdots f_{\zeta_{k+1}}f_{\zeta_{1}} \cdots f_{\zeta_k} \pi_{\lambda_2}$ is nonzero and $\lambda_1$-dominant.
\item[(2)] If $\lambda_2(\alpha_{\varepsilon_k}^\vee) \geq 1$, then $f_{\eta_1}\cdots f_{\eta_{q-1}}f_{\zeta_1}\cdots f_{\zeta_{r-1}} f_{\psi} f_{\varepsilon_p} \cdots f_{\varepsilon_{k+1}}f_{\varepsilon_{1}} \cdots f_{\varepsilon_k} \pi_{\lambda_2}$ is nonzero and $\lambda_1$-dominant.
\item[(3)] If $\lambda_2(\alpha_{\eta_k}^\vee) \geq 1$, then $f_{\varepsilon_1}\cdots f_{\varepsilon_p} f_{\zeta_1}\cdots f_{\zeta_{r-1}} f_{\psi} f_{\eta_{q-1}} \cdots f_{\eta_{k+1}} f_{\eta_1} \cdots f_{\eta_k} \pi_{\lambda_2}$ is nonzero and $\lambda_1$-dominant.
\item [(4)] If $\lambda_2(\alpha_{\psi}^\vee) \geq 1$, then $f_{\varepsilon_1} \cdots f_{\varepsilon_p}f_{\eta_1} \cdots f_{\eta_{q-1}} f_{\zeta_1} \cdots f_{\zeta_{r-1}} f_{\psi} \pi_{\lambda_2}$ is nonzero and $\lambda_1$-dominant.
\end{itemize}
\end{proof}

\begin{thm}\label{ThmSchurSupport}
Let $\lambda_1,\lambda_2,\lambda_3,\lambda_4$ be dominant weights of $\mathfrak{g}$ that satisfy hypothesis $(*)$. Let $(i_1,\cdots,i_t) \in I^t$ be a tuple satisfying the condition that $i_r\neq i_s$ if $r \neq s$ and the subdiagram of the Dynkin diagram of $\mathfrak{g}$ spanned by the nodes $i_1,\cdots,i_t$ is a disjoint union of connected components of type $A,D$, or $E$. Set $\beta=\alpha_{i_1}+\cdots+\alpha_{i_t}$ be an element of the root lattice. If $V(\lambda-\beta)$ is an irreducible component of $V(\lambda_3)\otimes V(\lambda_4)$ then $V(\lambda-\beta)$ is also an irreducible component of $V(\lambda_1)\otimes V(\lambda_2)$ i.e. $c_{\lambda_3,\lambda_4}^{\lambda-\beta} > 0$ implies that $c_{\lambda_1,\lambda_2}^{\lambda-\beta} > 0$.
\end{thm}

\begin{proof}
    Since, by Theorem \ref{thmSimplyLaced}, we may assume without loss of generality that the nodes $i_1,\cdots,i_t$ form a connected subdiagram of the Dynkin diagram of $\mathfrak{g}$, the proof follows immediately from proposition \ref{Schur} .
    
\end{proof}
\begin{rem}
The above theorem provides a new proof of a weaker form of the Schur positivity conjecture, originally established in \cite{DP} for type \(A\), in the special case when \(\nu = \lambda + \mu - \beta\), where \(\beta\) is a sum of distinct simple roots. In fact, this theorem extends the result to all simply laced finite simple Lie algebras, thereby offering further evidence in support of the Schur positivity conjecture.
\end{rem}

\begin{prop}

Let $\lambda_1,\lambda_2,\lambda_3,\lambda_4$ be four dominant integral weights of a simply laced symmetrizable Kac-Moody algebra $\mathfrak{g}$ such that $\lambda_1+\lambda_2=\lambda_3+\lambda_4=\lambda$. If for $1 \leq r \leq t$ and for simple roots $\alpha_{i_r}$, $k_r \leq \min\{\lambda_3(\alpha_{i_r}^\vee), \lambda_4(\alpha_{i_r}^\vee)\}$, then
$V(\lambda-\sum_{r=1}^t k_r \alpha_{i_r})$ is an irreducible component of $V(\lambda_3) \otimes V(\lambda_4)$. Moreover, if  
\[
\min\{\lambda_3(\alpha_j^\vee), \lambda_4(\alpha_j^\vee)\} \leq \min\{\lambda_1(\alpha_j^\vee), \lambda_2(\alpha_j^\vee)\}
\]
for all simple roots $\alpha_j$ and if $i_p \neq i_q$ for $p \neq q$, then $V(\lambda-\sum_{r=1}^t k_r \alpha_{i_r})$ is also an irreducible component of $V(\lambda_1) \otimes V(\lambda_2)$.

% That is, under the stated conditions on $\lambda_1, \lambda_2, \lambda_3, \lambda_4$ and the coefficients $k_r$, the existence of
% an irreducible component $V(\lambda - \sum_{r=1}^t k_r \alpha_{i_r})$ in $V(\lambda_3) \otimes V(\lambda_4)$ implies its existence
% in $V(\lambda_1) \otimes V(\lambda_2)$.

In particular, if $V(\lambda - k \alpha_i)$ is an irreducible component of $V(\lambda_3) \otimes V(\lambda_4)$ for some simple root $\alpha_i$ and $i \in I$, then $V(\lambda - k \alpha_i)$ is also an irreducible component of $V(\lambda_1) \otimes V(\lambda_2)$, and 
\[
c_{\lambda_3, \lambda_4}^{\lambda - k \alpha_i} \leq c_{\lambda_1, \lambda_2}^{\lambda - k \alpha_i}.
\]

\end{prop}

\begin{proof}
Note that by Lemma \ref{lemEquivDom}, the condition 
\[
k_r \leq \min\{\lambda_3(\alpha_{i_r}^\vee), \lambda_4(\alpha_{i_r}^\vee)\}
\]
implies that $f_{i_1}^{k_1} \cdots f_{i_t}^{k_t} \pi_{\lambda_4}$ is nonzero and $\lambda_3$-dominant. Our hypothesis 
\[
\min\{\lambda_3(\alpha_j^\vee), \lambda_4(\alpha_j^\vee)\} \leq \min\{\lambda_1(\alpha_j^\vee), \lambda_2(\alpha_j^\vee)\}
\]
implies that 
\[
k_r \leq \min\{\lambda_1(\alpha_{i_r}^\vee), \lambda_2(\alpha_{i_r}^\vee)\}.
\]
Hence, by Lemma \ref{lemEquivDom}, $f_{i_1}^{k_1} \cdots f_{i_t}^{k_t} \pi_{\lambda_2}$ is nonzero and $\lambda_1$-dominant, which yields 
\[
c_{\lambda_1, \lambda_2}^{\lambda - \sum_{r=1}^t k_r \alpha_{i_r}} \neq 0.
\]

If $V(\lambda - k \alpha_i)$ is an irreducible component of $V(\lambda_3) \otimes V(\lambda_4)$, then $f_i^k \pi_{\lambda_4}$ is nonzero and $\lambda_3$-dominant, since it is the unique element in $B(\lambda_4)$ of weight $\lambda_4 - k \alpha_i$. Moreover, for $i \neq j$,
\[
e_j f_i^k \pi_{\lambda_4} = 0,
\]
because the weight of $e_j f_i^k \pi_{\lambda_4}$ is $\lambda_4 - k \alpha_i + \alpha_j$, which is not less than $\lambda_4$, and also,
\[
e_i^\ell f_i^k \pi_{\lambda_4} = 0 \text{ if and only if } \ell > k.
\]
Thus $f_i^k \pi_{\lambda_4}$ is nonzero and $\lambda_3$-dominant if and only if 
\[
k \leq \min\{\lambda_3(\alpha_i^\vee), \lambda_4(\alpha_i^\vee)\}.
\]
By the hypothesis,
\[
\min\{\lambda_3(\alpha_i^\vee), \lambda_4(\alpha_i^\vee)\} \leq \min\{\lambda_1(\alpha_i^\vee), \lambda_2(\alpha_i^\vee)\},
\]
so $f_i^k \pi_{\lambda_2}$ is also nonzero and $\lambda_1$-dominant. The crystal element $\pi_{\lambda_1} \otimes f_i^k \pi_{\lambda_2}$ has weight $\lambda - k \alpha_i$, completing the proof.

\end{proof}

\subsection{Applications to Wahl components}
As an immediate application of our main theorem, we derive another proof of the existence of root components in the tensor product for certain special cases. These were previously established in \cite{Kumar1}, \cite{JK} for finite and affine Lie algebras using sophisticated techniques for \(N=1\), and in \cite{BJK} for arbitrary \(N\) but only for \(\mathfrak{sl}_n\).

\begin{cor}
Let \(\lambda\) and \(\mu\) be two dominant weights of a simply laced simple Lie algebra \(\mathfrak{g}\), and let \(R^+\) denote the set of positive roots of \(\mathfrak{g}\). Let \(N \geq 1\), and for a dominant weight \(\eta\) of \(\mathfrak{g}\), define 
\[
S_\eta = \{ i \in I : \eta(\alpha_i^\vee) < N \}.
\]
For any positive root \(\beta\) which is a sum of distinct simple roots, define
\[
F_\beta = \{ i \in I : \beta - \alpha_i \notin R^+ \cup \{0\} \}.
\]
Suppose that \(\lambda, \mu,\) and \(\beta \in R^+\) satisfy the following conditions:
\[
(1) \quad S_\lambda \cup S_\mu \subseteq F_\beta,
\]
\[
(2) \quad \lambda + \mu - N \beta \text{ is dominant.}
\]
A triple \((\lambda, \mu, \beta)\) satisfying the above conditions is called a \emph{Wahl triple} in the literature.

Then \(V(\lambda + \mu - N \beta)\) is an irreducible component of \(V(\lambda) \otimes V(\mu)\). Moreover, if \(\beta\) is not a simple root, then
\[
c_{\lambda, \mu}^{\lambda + \mu - N \beta} \geq 2.
\]
\end{cor}

\begin{proof}
Suppose \(\beta = \sum_{j=1}^t \alpha_{i_j} \in R^+\) satisfies the conditions of the hypothesis. We claim there exists a permutation \(\sigma \in S_t\) such that
\[
f_{i_{\sigma(1)}}^N \cdots f_{i_{\sigma(t)}}^N \pi_\mu
\]
is nonzero and \(\lambda\)-dominant.

For type \(A\), this follows from the fact that if 
\[
\beta = \alpha_u + \alpha_{u+1} + \cdots + \alpha_v
\]
with \(u \leq v\), then
\[
f_u^N f_{u+1}^N \cdots f_v^N \pi_\mu
\]
is nonzero and \(\lambda\)-dominant by Lemma \ref{lemEquivDom}, since the hypothesis on \(\beta\) implies 
\[
\lambda(\alpha_u^\vee) \geq N, \quad \lambda(\alpha_v^\vee) \geq N, \quad \mu(\alpha_u^\vee) \geq N, \quad \mu(\alpha_v^\vee) \geq N.
\]

Moreover, if \(\beta\) is not a simple root, Proposition \ref{LemEpsilon} states that
\[
\epsilon_u(f_u^N f_{u+1}^N \cdots f_v^N \pi_\mu) = N \quad \text{and} \quad \epsilon_u(f_v^N f_{v-1}^N \cdots f_u^N \pi_\mu) = 0.
\]
Hence, 
\[
f_u^N f_{u+1}^N \cdots f_v^N \pi_\mu \quad \text{and} \quad f_v^N f_{v-1}^N \cdots f_u^N \pi_\mu
\]
are distinct nonzero \(\lambda\)-dominant elements in \(B(\mu)\), which proves 
\[
c_{\lambda,\mu}^{\lambda+\mu-N\beta} \geq 2.
\]

It is well-known that a sum of distinct simple roots of a simple Lie algebra \(\mathfrak{g}\) is a root if and only if the Dynkin subdiagram formed by those simple roots is connected \cite[Proposition 4]{VV}.

Similarly to the type \(A\) case, in types \(D\) and \(E\), if the positive root 
\[
\beta = \alpha_{i_1} + \cdots + \alpha_{i_t}
\]
satisfies the hypotheses of the corollary, then 
\[
f_{u_1}^N \cdots f_{u_t}^N \pi_\mu
\]
is nonzero and \(\lambda\)-dominant, where \((u_1, \ldots, u_t) = (i_{\sigma(1)}, \ldots, i_{\sigma(t)})\) for some \(\sigma \in S_t\).
\end{proof}

\begin{rem}
Under the hypotheses of the above corollary, by the same argument it follows that for an arbitrary simply laced symmetrizable Kac–Moody algebra \(\mathfrak{g}\), if \(\lambda\) and \(\mu\) are dominant weights of \(\mathfrak{g}\) and \(\beta\) is a sum of distinct simple roots forming a connected subdiagram of the Dynkin diagram of \(\mathfrak{g}\), then \((\lambda, \mu, \beta)\) is a Wahl triple and 
\[
V(\lambda + \mu - \beta)
\]
is an irreducible component of \(V(\lambda) \otimes V(\mu)\). This follows from the fact that if 
\[
\beta = \alpha_{i_1} + \cdots + \alpha_{i_t}
\]
with \(i_1 < \cdots < i_t\) and the nodes \(i_1, \ldots, i_t\) form a connected subdiagram of the Dynkin diagram of \(\mathfrak{g}\), then \(\beta\) is a root, but \(\beta - \alpha_{i_k}\) is not a root for \(1 < k < t\), whereas both \(\beta - \alpha_{i_1}\) and \(\beta - \alpha_{i_t}\) are roots of \(\mathfrak{g}\). It is also true that \(\beta = \alpha_{i_1} + \cdots + \alpha_{i_t}\) is a root if and only if the nodes \(i_1, \ldots, i_t\) form a connected subdiagram of the Dynkin diagram of \(\mathfrak{g}\).
\end{rem}

\subsection{Nonsimply laced cases}

In this section, we work in the setting of Section 5 of \cite{Kashiwara95_2}, where Kashiwara introduces the notion of a virtualization map—though not under that name. More precisely, let $\mathfrak{g}$ be a complex non-simply-laced Kac-Moody Lie algebra, and let $\tilde{\mathfrak{g}}$ be a complex simply-laced Kac-Moody Lie algebra. We denote their Cartan matrices by $(a_{i,j})_{i,j \in I}$ and $(\tilde{a}_{i,j})_{i,j \in \tilde{I}}$, respectively. Assume there exists an injective map $\psi: I \to \tilde{I}$ such that for any dominant integral weight $\lambda$ of $\mathfrak{g}$, there exists a crystal $\tilde{\mathcal{B}}$ and a virtual crystal $\mathcal{V} \subset \tilde{\mathcal{B}}$ for which the crystal $B(\lambda)$ is isomorphic to $\mathcal{V}$ via a virtualization map
\[
\Upsilon_\lambda : B(\lambda) \to \mathcal{V}
\]
compatible with the folding $\psi$.

More concretely, if the $\tilde{\mathfrak{g}}$-crystal structure is given by 
\[
(\tilde{e}_i, \tilde{f}_i, \tilde{\epsilon}_i, \tilde{\phi}_i, \tilde{\mathrm{wt}})_{i \in \tilde{I}},
\]
then for each $i \in I$, the virtual crystal operators are defined by
\[
e_i^v = \prod_{j \in \psi(i)} (\tilde{e}_j)^{\gamma_i}, \quad f_i^v = \prod_{j \in \psi(i)} (\tilde{f}_j)^{\gamma_i}, \quad \epsilon_i^v = \gamma_i^{-1} \tilde{\epsilon}_j, \quad \phi_i^v = \gamma_i^{-1} \tilde{\phi}_j,
\]
for any $j \in \psi(i)$. This endows $\mathcal{V}$ with a virtual crystal structure isomorphic to the classical $\mathfrak{g}$-crystal structure
\[
(e_i, f_i, \epsilon_i, \phi_i, \mathrm{wt})_{i \in I}
\]
on $B(\lambda)$. The positive integers $(\gamma_i)_{i \in I}$ are called the \emph{scaling factors}; see \cite{AGHT} or \cite{BS} for their explicit values. Note also that via the virtualization map, the highest weight element $\pi_\lambda \in B(\lambda)$ corresponds to the highest weight element $\tilde{\pi}_\lambda \in \mathcal{V}$.

We present this virtualization framework following the recent preprint \cite{AGHT}. Furthermore, in \cite{PanScrimshaw}, the authors construct such maps explicitly using Littelmann's path model and provide necessary and sufficient conditions for a map to be a virtualization map. The virtualization procedure is well understood for finite-dimensional Lie algebras; see, for example, the book \cite{BS}. However, Kashiwara's article \cite{Kashiwara95_2} contains several examples in the infinite-dimensional setting.

\medskip
Thanks to the virtualization map, we can now prove an analogue of Proposition \ref{LemEpsilon}.

\begin{lem}\label{lemEpsilonNonSimply}
Let $\lambda$ be a dominant weight of $\mathfrak{g}$. Let 
\[
x = f_{i_1}^{b_1}\cdots f_{i_t}^{b_t} \pi_\lambda,
\]
where $i_p \neq i_q$ for $p \neq q$. Then, for $1 \leq r \leq t$, we have 
\[
\epsilon_{i_r}(x) = \gamma_{i_r}^{-1} \tilde{\epsilon}_k(\Upsilon_\lambda(x)),
\]
for any $k \in \psi(i_r)$. 

Moreover, if we denote the monomial 
\[
(f_{i_1}^v)^{b_1}\cdots (f_{i_t}^v)^{b_t} = \bigg(\prod_{j \in \psi(i_1)} (\tilde{f}_j)^{\gamma_{i_1}}\bigg)^{b_1} \cdots \bigg(\prod_{j \in \psi(i_t)} (\tilde{f}_j)^{\gamma_{i_t}}\bigg)^{b_t}
\]
by 
\[
(\tilde{f}_{j_1})^{c_1} \cdots (\tilde{f}_{j_s})^{c_s},
\]
where $c_u = b_u \gamma_{i_u}$, then
\[
\epsilon_{i_r}(x) = \gamma_{i_r}^{-1} \tilde{\epsilon}_k \big( (\tilde{f}_{j_1})^{c_1} \cdots (\tilde{f}_{j_s})^{c_s} \tilde{\pi}_\lambda \big) = \gamma_{i_r}^{-1} \max \{ 0, c_\ell + \sum_{u<\ell} c_u \tilde{a}_{j_u, j_\ell} \},
\]
for $1 \leq \ell \leq s$ such that $k = j_\ell \in \psi(i_r)$.
\end{lem}

\begin{proof}
The virtualization map $\Upsilon_\lambda$ is an isomorphism of crystals, so 
\[
\epsilon_{i_r}(x) = \epsilon_{i_r}^v(\Upsilon_\lambda(x)) = \gamma_{i_r}^{-1} \tilde{\epsilon}_k(\Upsilon_\lambda(x)),
\]
for any $k \in \psi(i_r)$.

Now,
\[
\Upsilon_\lambda(x) = (f_{i_1}^v)^{b_1} \cdots (f_{i_t}^v)^{b_t} \tilde{\pi}_\lambda = (\tilde{f}_{j_1})^{c_1} \cdots (\tilde{f}_{j_s})^{c_s} \tilde{\pi}_\lambda.
\]
Since the $i_p$ are distinct and so are the $j \in \psi(i_p)$, we may apply Proposition \ref{LemEpsilon} to obtain the stated formula.
\end{proof}

Now, Lemma \ref{lemAnyG} and Lemma \ref{lemEpsilonNonSimply} yield the following theorem.

\begin{thm}\label{thmNonSimplyLaced}
Let $t$ be a positive integer. Let $(b_1, \ldots, b_t) \in \mathbb{N}^t$ and $(c_1, \ldots, c_s) \in \mathbb{N}^s$ be tuples of integers as defined in Lemma \ref{lemEpsilonNonSimply}. Let $\lambda$ and $\mu$ be dominant integral weights of $\mathfrak{g}$ satisfying
\[
\lambda(\alpha_{i_t}^\vee) \geq b_t, \quad \mu(\alpha_{i_1}^\vee) \geq b_1,
\]
\[
\mu(\alpha_{i_r}^\vee) \geq \gamma_{i_r}^{-1} c_\ell + \sum_{u < \ell} \gamma_{i_r}^{-1} c_u \tilde{a}_{j_u, j_\ell}, \quad \text{for all } 1 < \ell \leq s, \ j_\ell \in \psi(i_r),
\]
\[
\lambda(\alpha_{i_r}^\vee) \geq b_r + \sum_{s > r} b_s a_{i_s, i_r}, \quad \text{for all } 1 \leq r < t.
\]
Then $V\big(\lambda + \mu - \sum_{p=1}^t b_p \alpha_{i_p}\big)$ is an irreducible component of $V(\lambda) \otimes V(\mu)$, provided $p \neq q$ implies $\alpha_{i_p} \neq \alpha_{i_q}$.
\end{thm}

\begin{proof}
The proof is analogous to the simply-laced case and is therefore omitted.
\end{proof}

\begin{rem}
Analogous applications to the three questions posed in the introduction can also be obtained in the nonsimply-laced case, but one must take into account the values of the scaling factors $(\gamma_i)_{i \in I}$.
\end{rem}

\section{Schur positivity conjecture for a simple Lie algebra $\mathfrak{g}$ in the case $\lambda \gg \mu$}
\label{se:Deep}

In this section, we assume that $\mathfrak{g}$ is a simple Lie algebra of rank $n$. Let $P$ be the weight lattice of $\mathfrak{g}$ and $P^+$ the set of dominant weights. We denote by $C_f$ the fundamental Weyl chamber of $\mathfrak{g}$. The Weyl group $W$ associated with $\mathfrak{g}$ acts naturally on $P$ and, by extension, on $P \otimes \mathbb{R}$.

We consider piecewise linear continuous paths, up to reparameterization,
\[
\pi : [0,1] \to P \otimes \mathbb{R}.
\]
For any choice of $\lambda \in \overline{C_f}$, $\pi_0 \in P \otimes \mathbb{R}$, $r \in \mathbb{N} \setminus \{0\}$, sequences $\underline{\tau} = (\tau_1, \tau_2, \dots, \tau_r)$ of elements in $W / W_\lambda$, and $\underline{a} = (a_0=0 < a_1 < \cdots < a_r = 1)$ in $\mathbb{R}$, we define a \emph{$\lambda$-path}, or a path of shape $\lambda$, $\pi = \pi(\lambda, \pi_0, \underline{\tau}, \underline{a})$ by
\[
\pi(t) = \pi_0 + \sum_{i=1}^{j-1} (a_i - a_{i-1}) \tau_i(\lambda) + (t - a_{j-1}) \tau_j(\lambda) \quad \text{for} \quad a_{j-1} \leq t \leq a_j.
\]
Any $\lambda$-path can be expressed in this way. We always assume that $\tau_j \neq \tau_{j+1}$.

Littelmann introduced a distinguished subset of $\lambda$-paths called \emph{LS paths} (after Lakshmibai and Seshadri), see \cite{Littelmann95}. Moreover, if $\lambda \in P^+$ is a dominant integral weight, the set of LS paths of shape $\lambda$ admits the structure of a crystal with crystal operators $e_i$ and $f_i$ defined geometrically. This crystal is isomorphic to the crystal $B(\lambda)$ associated with the highest weight irreducible representation $V(\lambda)$. The straight line path from $0$ to $\lambda$, denoted by $\pi_\lambda$, corresponds to the highest weight element in $B(\lambda)$. Elements in a tensor product of crystals $B_1 \otimes B_2$ correspond to concatenations of paths, denoted by $\pi * \eta$.

Furthermore, Littelmann showed that the set of LS paths of shape $\lambda$ can be characterized as the set of paths obtained from $\pi_\lambda$ by applying crystal lowering operators:
\[
f_{i_1} \cdots f_{i_m} \pi_\lambda.
\]
We will adopt this characterization as the definition of LS paths of shape $\lambda$.

\medskip

Let $\nu$ be a dominant weight. In the real vector space $V = X \otimes \mathbb{R}$, the convex hull of the Weyl group orbit $W \nu$ is a polytope called the \emph{Weyl polytope} of $\nu$, denoted by
\[
P_\nu^W = \operatorname{Conv} \{ w \nu \mid w \in W \}.
\]
Let $w \in W$ and $\beta$ be a simple root such that $\ell(s_\beta w) = \ell(w) + 1$. 
Then the segment $[w\nu, s_\beta w\nu]$ is either an edge of $P_\nu^W$ if the two points are distinct, or a single vertex if  $s_\beta w\nu = w\nu$. This edge can be labeled by the root $\beta$. The integer $w \nu(\beta^\vee)$ is called the \emph{$\beta$-length} of the segment, as in the final part of [\cite{Kam10}, Section 2.3].

\medskip

We begin by stating a known result to experts; the first author learned it from Travis Scrimshaw, while the second author thanks Peter Littelmann for pointing it out.

\begin{lem}
\label{lem:height_function}
Let $\lambda$ and $\lambda'$ be dominant weights such that $\lambda = \lambda' + \mu$, where $\mu$ is also dominant. Suppose the LS path
\[
\pi' = f_{i_t} \cdots f_{i_1} \pi_{\lambda'}
\]
exists. Then the LS path
\[
\pi = f_{i_t} \cdots f_{i_1} \pi_\lambda
\]
also exists. Moreover, if
\[
f_{i_t} \cdots f_{i_1} \pi_{\lambda'} \neq f_{j_t} \cdots f_{j_1} \pi_{\lambda'},
\]
then
\[
f_{i_t} \cdots f_{i_1} \pi_\lambda \neq f_{j_t} \cdots f_{j_1} \pi_\lambda.
\]
\end{lem}

\begin{proof}
Since all weights are dominant, the straight line path $\pi_\lambda$ and the concatenation $\pi_{\lambda'} * \pi_\mu$ generate isomorphic crystals. The tensor product rule for the operators $f_{i_j}$ implies they act only on the first factor:
\[
f_{i_t} \cdots f_{i_1} (\pi_{\lambda'} * \pi_\mu) = f_{i_t} \cdots f_{i_1} \pi_{\lambda'} * \pi_\mu.
\]
The lemma follows immediately.
\end{proof}

%\midskip
For any $\mu$ and $\nu$ two dominant weights such that $\mu + \nu = \lambda$. Consider 
$$
P_{\mu,\nu} = Conv\{\mu + w(\nu), w\in W\}.
$$ We have that $P_{\mu,\nu} \subseteq (\lambda -Q^+)\cap (\mu + w_0(\nu) + Q^+)$, with $Q^+ = \oplus_i \mathbb N\alpha_i$ and $w_0$ the longest element in the Weyl group. Further, one has $P_{\mu,\nu} = \mu + P^W_\nu$.

\begin{lem}
$P_{\lambda,\mu} \subset C_f$ if, and only if, $\lambda+w\mu \in C_f$ for all $w \in W$.
\end{lem}
\begin{proof}
By definition, $P_{\lambda,\mu} \subset C_f$ implies $\lambda+w\mu \in C_f$ for all $w \in W$. Let us assume that $\lambda+w\mu \in C_f$ for all $w \in W$. Then it is enough to prove that $\lambda+\sum_{i=1}^r a_i w_i\mu \in C_f$ for $a_i \in \mathbb{Q}^+$ and $w_i\in W$ for all $i$ satisfying the condition that $\sum_{i=1}^r a_i =1$. Let $w_k\mu(\alpha_j^{\vee})=\text{min}\{w_i\mu(\alpha_j^{\vee}): 1 \leq i \leq r\}$. Now $(\lambda+\sum_{i=1}^r a_i w_i\mu)(\alpha_j^{\vee}) \geq \lambda(\alpha_j^{\vee})+(\sum_{i=1}^r a_i)w_k\mu(\alpha_j^{\vee}) \geq (\lambda+w_k\mu)(\alpha_j^{\vee})\geq 0$.

\end{proof}

\begin{prop}
Let $\lambda_1, \lambda_2, \lambda_3, \lambda_4$ be dominant weights such that 
$\lambda_1 + \lambda_2 = \lambda = \lambda_3 + \lambda_4$.
Suppose that $\lambda_1$ and $\lambda_3$ are dominant enough in $C_f$ such that both polytopes $P_{\lambda_1,\lambda_2}$ and $P_{\lambda_3,\lambda_4}$ are contained in $C_f$. Then 
$$
V(\lambda_3)\otimes V(\lambda_4)\hookrightarrow V(\lambda_1)\otimes V(\lambda_2) \Longleftrightarrow P_{\lambda_3,\lambda_4}\subseteq P_{\lambda_1,\lambda_2}.
$$ 
\end{prop}

\begin{proof}
Assume that $V(\lambda_3)\otimes V(\lambda_4)\hookrightarrow V(\lambda_1)\otimes V(\lambda_2) $. Take a dominant weight $\nu\in P_{\lambda_3,\lambda_4}$ such that $\nu = \lambda_3 + \nu'$, where $\nu'\in P_{\lambda_4}^W \cap P$. Then since every element of $P_{\lambda_4}^{W} \cap P$ is a weight of $V(\lambda_4)$, there exists a LS path $\eta$ of shape $\lambda_4$ such that $\nu' = \eta (1)$. 

Since $P_{\lambda_3,\lambda_4}\subseteq C_f$, we get that $\nu = (\pi_{\lambda_3}*\eta) (1)$ is a dominant highest weight of $V(\lambda_3)\otimes V(\lambda_4)$. Then by our hypothesis, $\nu$ is also a dominant highest weight of $V(\lambda_1)\otimes V(\lambda_2)$. Therefore, there exists $\theta$, a LS path of shape $\lambda_2$ such that $\nu = (\pi_{\lambda_1}*\theta)(1)$. But $\theta(1) \in P^W_{\lambda_2}$ implies that $\nu \in P_{\lambda_1,\lambda_2}$. So we have $P_{\lambda_3,\lambda_4}\subseteq P_{\lambda_1,\lambda_2}$.

\medskip
Conversely, assume that $P_{\lambda_3,\lambda_4}\subseteq P_{\lambda_1,\lambda_2}$. Let $\mu = (\pi_{\lambda_3} * \theta) (1)$ be a dominant highest weight of $V(\lambda_3)\otimes V(\lambda_4)$, with $\theta$ a LS-path of shape $\lambda_4$. The weight $\mu$ is in $P_{\lambda_3,\lambda_4}$ and so $\mu\in P_{\lambda_1,\lambda_2}$ and hence $\mu$ is a dominant highest weight of the tensor product $V(\lambda_1)\otimes V(\lambda_2)$.

First we show that the conditions $P_{\lambda_3,\lambda_4}\subseteq P_{\lambda_1,\lambda_2}$ , $P_{\lambda_1,\lambda_2}\subseteq C_f$ and $P_{\lambda_3,\lambda_4}\subseteq C_f$ imply that $\lambda_2 \geq \lambda_4$. Indeed, from the containment in the fundamental chamber we get that $\lambda_1 - \lambda_2$ and $\lambda_3 - \lambda_4$ are dominant weights since for any simple reflection $s_i$, $\lambda_1 + s_i (\lambda_2)$ is dominant, so $(\lambda_1 + s_i (\lambda_2))(\alpha_i^\vee )\geq 0$. But  
$$
(\lambda_1 + s_i\lambda_2)(\alpha_i^\vee) = \lambda_1(\alpha_i^\vee) + ( s_i\lambda_2) (\alpha_i^\vee) = \lambda_1(\alpha_i^\vee) + \lambda_2(s_i\alpha_i^\vee) = (\lambda_1 - \lambda_2)(\alpha_i^\vee).
$$ Second, $P_{\lambda_3,\lambda_4}\subseteq P_{\lambda_1,\lambda_2}$ implies that the $\alpha_i$-length of the edge $[\lambda, \lambda_1 + s_i\lambda_2]$ is at least the $\alpha_i$-length of the edge $[\lambda, \lambda_3 + s_i\lambda_4]$. On the other hand, those lengths are equal to $(\lambda_1 + \lambda_2)(\alpha_i^\vee) - (\lambda_1 + s_i(\lambda_2))(\alpha_i^\vee)$ and $(\lambda_3 + \lambda_4)(\alpha_i^\vee) - (\lambda_3 + s_i(\lambda_4))(\alpha_i^\vee)$ respectively. Since $\lambda_1+\lambda_2=\lambda_3+\lambda_4$, we get that $(\lambda_3 - \lambda_4)(\alpha_i^\vee)\geq (\lambda_1 - \lambda_2)(\alpha_i^\vee)$. Hence $\lambda_3-\lambda_4 \geq \lambda_1 - \lambda_2$. Moreover $(\lambda_1 - \lambda_2 + 2\lambda_2) (\alpha_i^\vee) = \lambda (\alpha_i^\vee) = (\lambda_3 - \lambda_4 + 2\lambda_4) (\alpha_i^\vee)$ and $\lambda_3-\lambda_4 \geq \lambda_1 - \lambda_2$ says $\lambda_2\geq\lambda_4$.

\medskip

Now $\lambda_2$ and $\lambda_4$ are two dominant weights and the inequality $\lambda_2 \geq \lambda_4$ implies that $\lambda_2 = \lambda_4 + \mu$, for some dominant weight $\mu$.
So we are in the situation of applying Lemma \ref{lem:height_function} and if $\theta' = f_{i_t}\cdots f_{i_1} \pi_{\lambda_4}$, is a LS path of shape $\lambda_4$ then the LS path $f_{i_t}\cdots f_{i_1} \pi_{\lambda_2}$ exists and non zero. Moreover $f_{i_t}\cdots f_{i_1} \pi_{\lambda_4} \neq f_{j_t}\cdots f_{j_1} \pi_{\lambda_4}$ implies that $f_{i_t}\cdots f_{i_1} \pi_{\lambda_2} \neq f_{j_t}\cdots f_{j_1} \pi_{\lambda_2}$.  Therefore the multiplicity of the weight $\mu$ in $V(\lambda_3)\otimes V(\lambda_4)$ is less than or equal to the multiplicity of the weight $\mu$ in $V(\lambda_1)\otimes V(\lambda_2)$ i.e. $c_{\lambda_3,\lambda_4}^{\mu} \leq c_{\lambda_1,\lambda_2}^{\mu}$. 

\end{proof}

\begin{thm}
    Let $\lambda_1, \lambda_2, \lambda_3, \lambda_4$ be dominant weights such that 
$$
(*)
\left\{
\begin{array}{ll}
\lambda_1 + \lambda_2 = \lambda = \lambda_3 + \lambda_4  \\
\forall\alpha >0, \vert(\lambda_1-\lambda_2)(\alpha^\vee)\vert \leq \vert(\lambda_3 - \lambda_4)(\alpha^\vee)\vert .\\
\end{array}
\right.
$$ Suppose, in addition, that $\lambda_1$ and $\lambda_3$ are dominant enough inside the fundamental chamber $C_f$ such that $P_{\lambda_1,\lambda_2}$ and $P_{\lambda_3,\lambda_4}$ are contained in $C_f$. Then $P_{\lambda_3,\lambda_4}\subseteq P_{\lambda_1,\lambda_2}$. 
\end{thm}
\begin{proof}
As before, $P_{\lambda_1,\lambda_2} \subseteq C_f, P_{\lambda_3,\lambda_4} \subseteq C_f$ implies that $\lambda_1-\lambda_2$ and $\lambda_3-\lambda_4$ are two dominant weights. So the second condition in hypothesis $(*)$ reads: for any positive root $\alpha$, 
$$
0 \leq  (\lambda_1-\lambda_2)(\alpha^\vee) \leq (\lambda_3 -\lambda_4)(\alpha^\vee).
$$ 
Then, let us remark that, by convexity, it is enough to show that for any $w\in W$, $\lambda_3 + w\lambda_4$ sits inside $P_{\lambda_1,\lambda_2}$. If $\ell(w)=1$ then $\lambda_3+w\lambda_4=\lambda=\lambda_1+\lambda_2$ implies that $\lambda_3+w\lambda_4 \in P_{\lambda_1,\lambda_2}$.  Let $\beta$ be a simple root such that $v = s_\beta w$ be of length $\ell (w) + 1$. We prove by induction on the length of $w$ that the edge $[\lambda_3 + w\lambda_4, \lambda_3 + s_\beta w \lambda_4]$ is inside $P_{\lambda_1,\lambda_2}$. 

So, let $v = s_\beta$ for a simple root $\beta$. The edge starting from $\lambda$ labeled by $\beta$ connects the vertices $\lambda$ to $\lambda_1 + s_\beta(\lambda_2) = \lambda_1 + \lambda_2 -\lambda_2(\beta^\vee)\beta = \lambda -\lambda_2(\beta^\vee)\beta$. And we have 
$$
(\lambda_3 - \lambda_4 + 2\lambda_4)(\beta^\vee) = \lambda(\beta^\vee) = (\lambda_1 - \lambda_2 + 2\lambda_2)(\beta^\vee).
$$ So together with $(\lambda_1-\lambda_2)(\beta^\vee) \leq (\lambda_3 -\lambda_4)(\beta^\vee)$, we get $\lambda_2(\beta^\vee) \geq \lambda_4(\beta^\vee)$. Hence, the length of the edges starting from $\lambda$ labeled by $\beta$ of $P_{\lambda_1,\lambda_2}$ are at least those of $P_{\lambda_3,\lambda_4}$. So, $\lambda_3 + s_\beta\lambda_4$ lies inside $P_{\lambda_1,\lambda_2}$.

\medskip
Now assume that the edge $[\lambda_3 + u\lambda_4,\lambda_3 + w\lambda_4]$, for $w = s_\gamma u$, lies inside $P_{\lambda_1,\lambda_2}$ and let $v = s_\beta w$ be of length $\ell (w) + 1$. So $\alpha = w \beta$ is a positive root. And we also know that $\lambda_1 + w\lambda_2$, $\lambda_1 + w\lambda_2 - (w\lambda_2)(\beta^\vee)\beta = \lambda_1 + s_\beta w \lambda_2$ and the edge between these two vertices is in $P_{\lambda_1,\lambda_2}$. So, now $\lambda_3 + w\lambda_4 - (w\lambda_4)(\beta^\vee)\beta = \lambda_3 + s_\beta w \lambda_4$ is at the end of the edge labeled by $\alpha = w \beta$ and starting from $\lambda_3 + w\lambda_4$. But we have 
$$
(\lambda_3 - \lambda_4 + 2\lambda_4)(\alpha^\vee ) = \lambda(\alpha^\vee ) = (\lambda_1 - \lambda_2 + 2\lambda_2)(\alpha^\vee ).
$$ 
So again together with $(\lambda_1-\lambda_2)(\alpha^\vee ) \leq (\lambda_3 -\lambda_4)(\alpha^\vee )$, we get $\lambda_2(\alpha^\vee ) \geq \lambda_4(\alpha^\vee )$. Hence, the length of the edge $[\lambda_3 + w\lambda_4, \lambda_3 + s_\beta w \lambda_4]$ is at most the length of the parallel edge $[\lambda_1 + w\lambda_2, \lambda_1 + s_\beta w \lambda_2]$. So, for any $w$, by induction, the vertex $\lambda_3 + s_\beta w \lambda_4$ is at the end of a polyline built with edges parallel to the edges of $P_{\lambda_1,\lambda_2}$, but that are of lengths at most the lengths of the corresponding parallel edges in $P_{\lambda_1,\lambda_2}$, so the vertex $\lambda_3 + s_\beta w \lambda_4$ lies in $P_{\lambda_1,\lambda_2}$.

\end{proof}

Combining the above two results, we now have proved the Schur positivity conjecture under the extra assumption that $\lambda_3 > > \lambda_4$ and $\lambda_1 > > \lambda_2$.
\begin{rem}
From the above proof, we get that if $P_{\lambda,\mu} \subset C_f$, then 
$$
V(\lambda)\otimes V(\mu)=\bigoplus_{\nu \in \Pi(\mu)}V(\lambda+\nu)^{\oplus m_{\mu}(\nu)}
$$ where $\Pi(\mu)$ is the set of weights of $V(\mu)$ and $m_{\mu}(\nu)$ is the dimension of the weight space $V(\mu)_{\nu}$. This result is known by some experts as Travis Scrimshaw pointed out to us. 

Hence, the above theorem gives some non-trivial inequalities between dimensions of weight spaces of $V(\lambda_2)$ and $V(\lambda_4)$ as a corollary under the hypothesis of the theorem.
\end{rem}
Fix dominant integral weights $\lambda$, $\mu$ and an element $w$ of the Weyl group $W$ of the simple Lie algebra $\mathfrak{g}$. Let $v_{w\mu}$ denote a non-zero vector in the one-dimensional weight space of weight $w\mu$ in $V(\mu)$. The cyclic $\mathfrak{g}$-submodule $\mathcal{U}(\mathfrak{g})(v_{\lambda}\otimes v_{w\mu})$ of $V(\lambda)\otimes V(\mu)$, where $\mathcal{\mathfrak{g}}$ denotes the universal enveloping algebra of $\mathfrak{g}$, is called a Kostant-Kumar module which has been studied in \cite{KRV}. 
From the above proof, the following corollary for Kostant-Kumar modules is obvious if we use the formulae given in \cite[Theorem 6.1]{KRV} for the decomposition of Kostant-Kumar modules into direct sum of irreducible modules.
\begin{cor}
    Let $\lambda_1, \lambda_2, \lambda_3, \lambda_4$ be dominant weights such that 
$$
(*)
\left\{
\begin{array}{ll}
\lambda_1 + \lambda_2 = \lambda = \lambda_3 + \lambda_4  \\
\forall\alpha >0, \vert(\lambda_1-\lambda_2)(\alpha^\vee)\vert \leq \vert(\lambda_3 - \lambda_4)(\alpha^\vee)\vert .\\
\end{array}
\right.
$$
Suppose that $\lambda_1$ and $\lambda_3$ are dominant enough in $C_f$ such that both polytopes $P_{\lambda_1,\lambda_2}$ and $P_{\lambda_3,\lambda_4}$ are contained in $C_f$.
Then for a Weyl group element $w$ of the simple Lie algebra $\mathfrak{g}$, there is an injective $\mathfrak{g}$-module homomorphism from $\mathcal{U}\mathfrak{g}(v_{\lambda_3}\otimes v_{w\lambda_4})$ to $\mathcal{U}\mathfrak{g}(v_{\lambda_1}\otimes v_{w\lambda_2})$ where $v_{w\lambda_i}$ is the extremal weight vector of weight $w\lambda_i$ inside $V(\lambda_i)$ and $\mathcal{U}\mathfrak{g}$ is the universal enveloping algebra of $\mathfrak{g}$.
\end{cor}
\begin{proof}
This just follows from Theorem 6.1 in \cite{KRV} and the observation that for $\{i_1,\cdots,i_t\} \subseteq I$ if $f_{i_1}\cdots f_{i_t}v_{w\lambda_4} \in B_w(\lambda_4)$, then $f_{i_1}\cdots f_{i_t}v_{w\lambda_2} \in B_w(\lambda_2)$ where $B_w(\lambda_i)$ is the Demazure subcrystal inside the crystal $B(\lambda_i)$ corresponding to the Weyl group element $w$ for $i=2,4$.

\end{proof}

\medskip
{\bf Acknowledgments. } The authors thank Peter Littelmann and Travis Scrimshaw for valuable discussions and Katsuyuki Naoi, Shrawan Kumar, Prakash Belkale, Arun Ram and Travis Scrimshaw for going through the paper very carefully and making comments on an earlier version. R.B. is grateful to Sankaran Viswanath for his patience in trying to understand the problem studied in this paper and giving insightful suggestions particularly to the proof of Proposition \ref{Schur} when she was visiting The Institute of Mathematical Sciences(IMSc) in Chennai in Jan 2025 and she also thanks IMSc for excellent working conditions at which part of this work was carried out. R.B. also thanks Amartya Shekhar Dubey for his help with the Dynkin diagram drawn in the paper.The authors thank the referee, whose comments helped to improve the presentation of the article.

\medskip
{\bf Author contributions. } All authors equally contributed to the manuscript.

\medskip
{\bf Funding. } The authors also thank the Université Jean Monnet in Saint-Étienne for the one month professorship position accorded to R.B when this collaboration was initiated.

\bibliographystyle{amsplain}

\end{document}